\documentclass[a4paper,11pt]{article}
\usepackage{amssymb}

\newcommand{\qeda}{\hspace{10mm}\hfill $\square$}

\newtheorem{theorem}{Theorem}
\newtheorem{lemma}{Lemma}
\newtheorem{proposition}{Proposition}
\newtheorem{remark}{Remark}

\title{Pseudo-radial solutions of semi-linear elliptic equations\\ on symmetric domains}
\author{Ahmad El Soufi and Mustapha Jazar\footnote{The first author is supported by a grant from the Lebanese University.}}

\date{}
\begin{document}

\maketitle

\begin{center}
{\bf\small Abstract} \vspace{3mm}

\hspace{.05in}\parbox{4.5in} {\small In this paper we investigate
existence and characterization of non-radial pseudo-radial (or
separable) solutions of some semi-linear elliptic equations on
symmetric 2-dimensional domains. The problem reduces to the phase
plane analysis of a dynamical system. In particular, we give a
full description of the set of pseudo-radial solutions of
equations of the form $\Delta u = \pm a^2(|x|) u|u|^{q-1}$, with
$q>0$, $q\neq 1$. We also study such equations over spherical or
hyperbolic symmetric domains.}
\end{center}



\section{Introduction}

A very rich literature has been devoted during the last decades to
the study of semi-linear elliptic equations of the form

\begin{equation}\label{edpq}
\Delta u=\varepsilon u|u|^{q-1}
\end{equation}
over a symmetric Euclidean domain $\Omega$, where $\varepsilon=\pm 1$ and $q$ is a
positive real number, $q\ne 1$.
 In particular, it is known that when the domain $\Omega$ is an Euclidean ball,
then any positive solution of (\ref{edpq}) with $q>1$, satisfying
Dirichlet boundary conditions, is radial (see \cite{GNN}, see also
\cite{GNN1} for the case $\Omega={\mathbb R}^n$), while on an
annular domain $\Omega=\left\{0<R<|x|<R+ c \right\}$, the equation
(\ref{edpq}) admits positive non-radial solutions for $\varepsilon =-1$ and
any $q>1$ (see \cite{BN,C,NS1, S}). See also the important work of Rabinowitz \cite{R} where sufficient conditions are given for the existence of infinitely many non necessarily radial solutions to such equations.

The existence of non-radial solutions is generally based on
minimization techniques (see \cite{S} for a survey). Another natural approach toward the
existence of non-radial solutions consists in searching solutions
of the form
$$u(r,\theta)=h(r)w(\theta),$$
where $(r,\theta)$ are polar coordinates and $w$ is a
$2\pi$-periodic function. Such a solution is sometimes called {\it
``separable"} or {\it ``pseudo-radial''} (see for instance
\cite{KV}). Of course, $u$ is non-radial as far as $w$ is
nonconstant. Moreover, we say that $u$ is of {\it mode $k
\in{\mathbb N} $} if the least period of $w$ is $\frac{2\pi}{k}$
(Thus, radial solutions could be considered as being of mode
$+\infty$). Notice that singular pseudo-radial solutions
play an important role in the study and classification of singularities of
solutions of semi-linear equations.

In the present work, we make use of this last approach in order to
investigate the existence of possibly singular non-radial solutions to
\begin{itemize}
\item[(a)] ``weighted'' equations of the form
\begin{equation}\label{conf}
\Delta u= \varepsilon a^2(r)u|u|^{q-1},
\end{equation}
on a rotationally symmetric (non necessarily bounded) Euclidean domain $\Omega\subset{\mathbb R}^2$
(see for instance \cite{Na,
Ni, NiN} for results concerning radial solutions),

\item[(b)] equation (\ref{edpq}) over
rotationally symmetric non necessarily flat domains like geodesic
discs or annular domains of the sphere ${\mathbb S}^2$, the
hyperbolic plane ${\mathbb H}^2$, a revolution surface, etc.

\end{itemize}
Notice that our study includes the ``sublinear'' case $0<q<1$
which is of particular interest since it is not much treated in
the literature.

The natural general setting in which these topics can be put is
the following. We consider Equation (\ref{edpq}) on $I\times
\mathbb{S}^{1}$, where $I$ is an interval, $I\subset (0,\infty)$,
endowed with a Riemannian metric of the form
$$g= a^2(r)dr^2+b^2(r)d\theta^2,$$
where $a$ and $b$ are two positive differentiable functions on $I$.
The Laplace operator associated with the metric $g$ is given by
\begin{equation}\label{Deltag}
\Delta^g=\frac 1{a^2}\frac {\partial^2}{\partial r^2}+\frac
1{ab}\left(\frac{b}{a}\right)'\frac\partial{\partial r}+ \frac 1{b^2}\frac
{\partial^2}{\partial \theta^2}.
\end{equation}

Recall that a rotationally symmetric domain $\Omega\subset
{\mathbb R}^2$, ${\mathbb S}^2$ or ${\mathbb H}^2$ can be
identified, using appropriate polar coordinates $(r,\theta)$, with
a cylinder $I\times \mathbb{S}^{1}$ endowed with the Riemannian
metric $g=dr^2+b^2(r)d\theta^2$, with
$$b(r)=\left\{\begin{array}{ll} r & \mbox{in the Euclidean case}\\
\sin r&\mbox{in the spherical case}\\ \sinh r&\mbox{in the
hyperbolic case.}
\end{array} \right.$$
Here, the $r$ variable represents in each case the geodesic
distance to the center of the domain. On the other hand, the
weighted equation $\Delta u= \varepsilon a^2(r)u|u|^{q-1}$ of item
(a) above is equivalent to the equation
$$\Delta^gu=\varepsilon u|u|^{q-1},$$
where $g= a^2(r) g_{{\mathbb R}^2}=a^2(r) (dr^2+r^2d\theta^2)$ is a
Riemannian metric conformal to the Euclidean one $g_{{\mathbb R}^2}$.

Hence, we consider for any Riemannian metric $g=
a^2(r)dr^2+b^2(r)d\theta^2$ on $I\times \mathbb{S}^{1}$, the PDE
\begin{equation}\label{edp}
\Delta^gu=\varepsilon u|u|^{q-1},
\end{equation}
where $\varepsilon=\pm1$, and look for non-radial pseudo-radial
solutions. We first show (Theorem \ref{main}) that a necessary and
sufficient condition for equation (\ref{edp}) to admit such a
non-radial pseudo-radial solution consists in the existence of a
real number $\mu$ such that:
\begin{itemize}
\item[(i)] (condition on the metric $g$)
\begin{equation}
{\left( a^{-1}b^{\frac 2{1-q}}b'\right)}'=\mu\frac{(1-q)}2 a
b^{\frac {1+q}{1-q}},
\end{equation}
\item[(ii)] the following ODE admits a nonconstant $2\pi$-periodic
solution
\begin{equation}\label{ode2}
w''(\theta)+\mu w(\theta)=\varepsilon w(\theta)|w(\theta)|^{q-1}.
\end{equation}
\end{itemize}
Moreover, when (i) and (ii) are satisfied, then
$u(r,\theta)=b^{\frac 2{1-q}}(r)w(\theta)$ is a solution of
(\ref{edp}). Notice that a more general version is actually given
in Theorem \ref{main}.

For example, condition (i) is satisfied for the Euclidean metric
($a=1$ and $b=r$) with $\mu=\frac 4{(1-q)^2}$, the spherical
metric ($a=1$ and $b=\sin r$) with $\mu=1$ and $q=3$, and the
hyperbolic metric ($a=1$ and $b=\sinh r$) with $\mu=1$ and $q=3$.
On the other hand, in the conformal case (i.e. $b(r)=ra(r)$), a
consequence of condition (i) is that, if the weighted equation
$\Delta u=\varepsilon a^2(r)u|u|^{q-1}$ admits a non-radial
pseudo-radial solution in a rotationally symmetric $\Omega\subset {\mathbb R}^2$,
then the function $a(r)$ has one of the
following forms (Theorem \ref{th2}):
$$a(r)=r^{-1}(Mr^\alpha+Nr^{-\alpha})^{\frac{1-q}2},$$

$$a(r)=r^{-1}(M+N\ln r)^{\frac{1-q}2},$$
or
$$a(r)=r^{-1}{\left[ M\cos(\alpha\ln r) +
N\sin(\alpha\ln r) \right]}^{\frac{1-q}2}$$
with $M$, $N \in{\mathbb R}$ and $\alpha>0$.

Condition (ii) leads us to study the ODE (\ref{ode2}) and seek its
$2\pi$-periodic solutions, according to
the values of $\varepsilon$, $q$ and $\mu$. Actually, this study
constitutes the main part of this paper. The case
$\varepsilon=+1$, $0<q<1$ and $\mu>0$ is the more interesting and
novel one because of the lack of regularity of the nonlinear term
at the origin. Notice that some earlier results concerning the ODE
(\ref{ode2}) have been obtained for particular values of
$\varepsilon$, $q$ and $\mu$ (see for instance \cite{BVB, CI, Ve1}).

In particular, for the weighted equation
\begin{equation}\label{weighted1}
\Delta u=\varepsilon r^{-2}(Mr^\alpha+Nr^{-\alpha})^{1-q}u|u|^{q-1},
\end{equation}
in ${\mathbb R}^2\backslash\left\{O\right\}$,
the corresponding $\mu$ is equal to $\alpha^2$ and we obtain (Theorem \ref{th3}), as an application
of the ODE analysis, that~:
\begin{itemize}
\item if $\varepsilon= -1$, then, for every $ q\ne1$, $\alpha>0$,
equation (\ref{weighted1}) admits for any integer $k>\alpha$, a
unique (up to sign) pseudo-radial solution of mode $k$ which is
sign changing,
 \item if $\varepsilon= +1$, $ q>1$ and $\alpha>1$, then for every $ k\ge1$, equation (\ref{weighted1})
admits a unique (up to sign) pseudo-radial solution of mode $k$ which is sign changing,
 \item if $\varepsilon= +1$, $ 0<q<1$ and $\alpha>1$, then equation (\ref{weighted1}) admits,
for any integer $k\in(\frac{1-q}{2}\alpha, \alpha)$,
 a unique (up to sign) sign changing pseudo-radial solution of mode $k$, and,
for any possible integer $k\in ((1-q)\alpha ,\sqrt{1-q}\, \alpha)$, a unique positive
pseudo-radial solution of mode $k$. Moreover, if $(1-q)\alpha$ is an integer, then there exists a
unique nonnegative pseudo-radial solution of mode $(1-q)\alpha$,
\item if $\varepsilon= +1$ and $0<\alpha\le1$, then any pseudo-radial solution of (\ref{weighted1}) is radial.
\end{itemize}
In all cases, the radial part of pseudo-radial solutions is given by $h(r)=Mr^\alpha+Nr^{-\alpha}$.

In a similar manner, we describe the set of
non-radial pseudo-radial solutions of the weighted equations (see Theorems \ref{th4} and \ref{th5})
$$\Delta u=\varepsilon r^{-2}{\left[ M+N\ln
r\right]}^{1-q}u|u|^{q-1}$$
and
$$\Delta u=\varepsilon r^{-2}{\left[ M\cos(\alpha\ln r) +
N\sin(\alpha\ln r)\right]}^{1-q}u|u|^{q-1}.$$

Concerning the equation $\Delta_{{\mathbb S}^2} u=\varepsilon
u|u|^{q-1}$ (resp. $\Delta_{{\mathbb H}^2} u=\varepsilon
u|u|^{q-1}$) on a rotationally symmetric domain of the sphere
(resp. the hyperbolic plane), it turns out that any pseudo-radial
solution is radial unless for $\varepsilon=+1$ and $q=3$ where,
for any integer $k\ge 2$, there exists a unique (up to sign)
pseudo-radial solution of mode $k$ which is sign changing.

The paper is organized as follows. In the second section we prove
Theorem \ref{main} giving necessary and sufficient conditions for
equation (\ref{edp}) to admit a non-radial pseudo-radial solution.
The third section is devoted to the study of the resulting ODE
(\ref{ode2}). In the last section we apply our results to answer
completely the question of existence of positive and sign-changing
non-radial pseudo-radial solutions in the conformal case as well
as the case of symmetric domains of standard spaces.

\section{A general result}

Let $\Omega\subset {\mathbb R}^N$ be a rotationally symmetric
domain (e.g. a ball, a spherical shell, ${\mathbb R}^N
\backslash\{O\}$, etc.) that we parametrize by spherical
coordinates $(r,\sigma)\in I\times {\mathcal{S}^{N-1}}\mapsto
r\sigma \in \Omega$. We endow $\Omega$ (or, equivalently, $I\times
{\mathcal{S}^{N-1}}$) with the Riemannian metric
$$g= a^2(r)dr^2+b^2(r)d\sigma^2,$$
where $a$ and $b$ are two positive differentiable functions on $I$.
The Laplace-Beltrami operator associated with the
metric $g$ is given by
\begin{equation}\label{Deltag1}
\Delta^g=\frac 1{b^2}\left[c^2\frac {\partial^2}{\partial
r^2}+(N-1)c c'\frac\partial{\partial r}+\Delta_{\mathcal{S}^{N-1}}
\right],\end{equation} where $c^2=\frac {b^2}{a^2}$ and
$\Delta_{\mathcal{S}^{N-1}}$ stands for the Laplace-Beltrami operator of
the standard sphere. Consider the equation
\begin{equation}\label{edpN}
\Delta^gu=\varepsilon u|u|^{q-1},
\end{equation}
where $\varepsilon=\pm 1$ and $q$ is a positive real number, $q\ne
1$. A solution $u$ in $\Omega$ is said to be ``pseudo-radial'' if
it can be written, with respect to the $(r,\sigma)$ coordinates,
as
$$u(r,\sigma)=h(r)w(\sigma).$$

\begin{theorem}\label{main}
Equation (\ref{edpN}) admits a non radial pseudo-radial solution if and
only if there exists $\mu\in{\mathbb R}$ such that:
\begin{itemize}
\item[(i)] (Condition on the metric $g$)
\begin{equation}\label{conditiongg}
c^2(b^{\frac 2{1-q}})''+ (N-1)cc'(b^{\frac 2{1-q}})'=\mu b^{\frac 2{1-q}}.
\end{equation}
\item[(ii)] the following equation admits a nonconstant solution
on $\mathbb{S}^{N-1}$
\begin{equation}\label{conditionS}
\Delta_{\mathbb{S}^{N-1}}w+\mu w=\varepsilon w|w|^{q-1}.
\end{equation}
\end{itemize}
Moreover, when (i) and (ii) are satisfied, then any non radial
pseudo-radial solution of (\ref{edpN}) is of the form
$u(r,\sigma)=b^{\frac 2{1-q}}(r)w(\sigma)$, where $w$ is a
nonconstant solution of (\ref{conditionS}).
\end{theorem}

Notice that in dimension N=2, the condition on the metric (\ref{conditiongg}) reads
\begin{equation}\label{conditiong}
{\left( a^{-1}b^{\frac 2{1-q}}b'\right)}'=\frac{\mu(1-q)}2 a
b^{\frac {1+q}{1-q}},
\end{equation}
while the equation (\ref{conditionS}) reduces to the ODE (\ref{ode2})
$$
w''(\theta)+\mu w(\theta)=\varepsilon w(\theta)|w(\theta)|^{q-1},
$$
where $w$ is a $2\pi$-periodic function.

\medskip

The proof of Theorem \ref{main} relies on the following elementary
lemma.

\begin{lemma}
Let $a$, $b$ and c be three nontrivial real-valued functions on a
set $X$ and let $\alpha$, $\beta$ and $\gamma$ be three
differentiable functions on an interval $I\subset {\mathbb R}$
such that $\gamma$ admits no zeros in $I$. Assume that, for every
$ (x, y)\in I\times X$,
\begin{equation}\label{prop}
\alpha(x)a(y)+\beta(x)b(y)=\gamma(x)c(y),
\end{equation}
then,
 \begin{itemize}
 \item either $\alpha$ and $ \beta$ are
proportional to $\gamma$ on $I$,
\item or
$a$, $b$ and $c$ are mutually proportional on $X$.
\end{itemize}
\end{lemma}
\noindent\textbf{Proof.} Dividing (\ref{prop}) by $\gamma (x)$
and, then, differentiating with respect to the $x$ variable, one
obtains, for every $ (x,y)\in I\times X$,
\begin{equation}\label{prop1}
\bar \alpha (x) a(y)+\bar\beta (x) b(y)=c(y)
\end{equation}
and
\begin{equation}\label{prop2}
\bar \alpha' (x) a(y)=-\bar\beta'(x) b(y),
\end{equation}
where $\bar\alpha:={\alpha\over\gamma}$ and $\bar\beta:={\beta\over\gamma}$. Two cases are to be considered.

\begin{enumerate}
\item Assume that $ab \equiv 0$. Multiplying (\ref{prop2}) by $a$
and, then, by $b$, we deduce that $ \bar \alpha' =\bar
\beta'\equiv 0$. Hence, $\bar\alpha$ and $\bar\beta$ are constants
which means that $\alpha$ and $\beta$ are proportional to $\gamma$
on $I$. \item Assume that there exists $y_0\in X$ such that
$a(y_0)b(y_0)\ne 0$. Setting $K=b(y_0)/a(y_0)$, we deduce from
(\ref{prop2}) that $\bar\alpha'(x)=-K\bar\beta'(x)$, that is
$\bar\alpha(x)=-K\bar\beta(x)+C$ for some $C\in {\mathbb R}$, and,
either $ \bar \alpha' =\bar \beta'\equiv 0$ or $b(y)=Ka(y)$.
Therefore, using (\ref{prop1}), either $\alpha$ and $\beta$ are
proportional to $\gamma$ on $I$, or $b$ and $c$ are proportional
to $a$ on $X$.
\end{enumerate}

$\hfill\square$

\noindent\textbf{Proof of Theorem \ref{main}.}
 If $u(r,\sigma)=h(r)w(\sigma)$ is a non-radial pseudo-radial solution of (\ref{edpN}), then
 \begin{eqnarray}\label{separate}
\left[c^2h''(r)+(N-1) c c'h'(r)\right] w(\sigma)
+h(r)\Delta_{\mathcal{S}^{N-1}}w(\sigma) &=&\hfill\nonumber\\
&&\hskip-6cm\varepsilon
b^2(r)h(r)|h(r)|^{q-1}w(\sigma)|w(\sigma)|^{q-1}.
\end{eqnarray}

Let us apply the last lemma on $I_1\times {\mathbb S}^{N-1}$,
where $I_1$ is a subinterval on which the function $h $ has no zeros.
 Since $w$ is not constant ($u$ is assumed to be non-radial) and $q\ne 1$,
the function $w|w|^{q-1}$ can never be proportional to $w$ on ${\mathbb S}^{N-1}$.
Hence, there exist two constants $\lambda$ and $\mu$ such that
\begin{equation}\label{eq1} h(r)=\lambda
b^2(r)h(r)|h(r)|^{q-1},\end{equation}
\begin{equation}\label{eq2} c^2h''+ (N-1)cc'h'=\mu h
,\end{equation} and, hence,
\begin{equation}\label{edpS}\Delta_{\mathcal{S}^{N-1}}w(\sigma)+\mu w(\sigma)=\frac
\varepsilon\lambda w(\sigma)|w(\sigma)|^{q-1}.\end{equation}

From (\ref{eq1}) we have $|h(r)|={(\lambda b^2)}^{-\frac 1{q-1}}$. Hence,
the function $|h|$ is proportional to ${b}^{-\frac 2{q-1}}$ as long as it does not vanish.
For continuity reasons, this implies that $|h(r)|={(\lambda b^2)}^{-\frac 1{q-1}}$ on the whole $I$.
In particular, $h$ does not vanish on $I$ and one can assume, without loss of generality, that
$$h={b}^{-\frac 2{q-1}}.$$
Indeed, it is clear that if $u$ is a solution of
(\ref{edp}) then $-u$ is also a solution. On the other hand,
replacing $h$ by
$\lambda^{-\frac 1{q-1}}h$ and $w$ by $\lambda^{\frac 1{q-1}}w$,
the solution $u$ remains unchanged and the PDE (\ref{edpS})
reduces to (\ref{conditionS}). Finally, replacing $h$ by ${b}^{-\frac 2{q-1}}$ in (\ref{eq2}) one gets (\ref{conditiong}).

Conversely, it is easy to check that if conditions (i) and (ii)
are satisfied, then the function $u(r,\sigma)=b^{\frac
2{1-q}}(r)w(\sigma)$, where $w$ is a nonconstant solution of
(\ref{conditionS}), is a solution of (\ref{edpN}). \qeda

\section{Study of the ODE}
In this section we investigate the existence of $2\pi$-periodic
solutions of the ODE (\ref{ode2})
$$
w''(\theta)+\mu w(\theta)=\varepsilon w(\theta)|w(\theta)|^{q-1}
$$
according to the values of the parameters $\varepsilon$, $q$ and $\mu$.

In order to transform the ODE into a dynamical system we put
$x=w$ and $y=w'$. We then get
$$(\mathcal{S})\left\{
\begin{array}{l}
\,\,x'=P(y):=y,\\
\\
y'=Q(x):=-\mu x+\varepsilon x|x|^{q-1}.
\end{array}
\right.$$
The origin is either a critical point ($q>1$) or the only singular point ($0<q<1$) of the system $(\mathcal{S})$.
Notice that solutions of (\ref{ode2}) satisfy
\begin{equation}\label{ip}
w'^2(t)=w'^2(0)-\mu\left(w^2(t)-w^2(0)\right)+\frac {2\varepsilon}{q+1}
\left[ |w(t)|^{q+1}-|w(0)|^{q+1} \right].
\end{equation}
Equivalently, the orbit of $(\mathcal{S})$ passing through the point $(x_0,y_0)$ is given by the equation
$$y^2-y_0^2=-\mu \left( x^2-x_0^2\right)+\frac{2\varepsilon}{q+1}\left[|x|^{q+1}-|x_0|^{q+1} \right].$$

If an orbit of $(\mathcal{S})$ intersects one of the coordinates
axes, then the intersection occurs perpendicularly (indeed, if
$(x(t),y(t))$ satisfies $(\mathcal{S})$, then $x(t_0)=0$ implies
$y'(t_0)=0$ and $y(t_0)=0$ implies $x'(t_0)=0$). Moreover, the
system is clearly invariant under the transformations
$\Phi_x:(t,x,y)\mapsto (-t,-x,y)$ and $\Phi_y:(t,x,y)\mapsto
(-t,x,-y)$. The following lemma is a quasi-immediate consequence
of these observations.

\begin{lemma}\label{l2}
The coordinates axes are axes of symmetry of the dynamical
system $(\mathcal{S})$ and the origin is a center of
symmetry.

Any orbit which intersects both the two coordinates axes is
necessarily closed.

\end{lemma}

\subsection{Case $q>1$ and $\varepsilon=-1$}

\subsubsection{Assume $\mu\ge0$}

A classical study of the phase plane of the system $(S)$ in this
case shows that all solutions are periodic that turn around the
origin. This can be seen considering for example the following
Lyapounov function
$$E(w,w'):=\mu w^2+\frac 2{q+1}|w|^{q+1}+w'^2.$$
Moreover, periodic solutions are given by the equation

$$E(x,y)=y_0^2$$
where $y_0=y(0)=w'(0)$.

\begin{lemma}\label{pro1}
The period function $y_0\in (0,\infty) \mapsto T(y_0)$, where
$T(y_0)$ is the period of the solution of
$(\mathcal{S})$ passing through the point $(0,y_0)$, is
decreasing with
$$T(0,\infty)=\left\{\begin{array}{ll} (0, \frac{2\pi}{\sqrt \mu}) & \mbox{if $\mu >0$}\\
(0,\infty)&\mbox{if $\mu = 0$.}
\end{array} \right.$$
\end{lemma}
\noindent{\bf Proof.} Notice that $P$ is homogeneous and
$Q$ is sub homogeneous of degree 1 (indeed, $Q(\nu\,x)-\nu Q(x)=(\nu-\nu^q)x^q<0$ for all $x>0$ and
all $\nu>1$). Applying \cite[Theorem 3]{HJV}, we get the monotony of
the period function.

Let $\mathcal{O}$ be the orbit of a periodic solution $w$ such
that $w(0)=0$ and $w'(0)=s$, with $s>0$. For symmetry reasons, it
suffices to work with the quarter of the orbit lying in the region
$(x\ge 0)\cap(y\ge 0)$. In this region, one has, from (\ref{ip}),
$w'(t)=\sqrt{s^2-U(w(t))}$, where $U(x)=\mu x^2+\frac
2{q+1}x^{q+1}$. Denoting by $(x(s),0)$ the first intersection
point of the orbit with the $x$-axis, we have $w(T(s)/4)=x(s)$,
$U(x(s))=s^2$ and, then,
\begin{eqnarray*}
T(s)&=&4\int_0^{x(s)}\frac{dx}{\sqrt{s^2-U(x)}}=4\int_0^1\frac{x(s)d\tau}{\sqrt{U(x(s))-U(x(s)\tau)}}\\
&=& 4\int_0^1\frac{d\tau}{\sqrt{(1-\tau^2)\mu+\frac
2{q+1}x(s)^{q-1}(1-\tau^{q+1})}}.\\
\end{eqnarray*}
Since $x(s)\to0$ as $s\to0$, and $\int_0^1\frac{d\tau}{\sqrt{(1-\tau^2)}} =
\frac{\pi}{2}$, we
deduce using standard convergence results,
$$\lim_{s\to 0}T(s)=\left\{\begin{array}{ll} \frac{2\pi}{\sqrt \mu} & \mbox{if $\mu >0$}\\
\infty&\mbox{if $\mu = 0$.}
\end{array} \right.$$
Since $x(s)\to \infty$ as $s\to \infty$, we obtain using the same calculations,
$$\lim_{s\to\infty}T(s)=0.$$ \qeda

A direct consequence of Lemma \ref{pro1} is the following
\bigskip

\begin{proposition}\label{pro2}
Assume that $q>1$, $\mu\ge0$ and $\varepsilon=-1$. Then for all
integers $k>\sqrt\mu$ the ODE (\ref{ode2}) admits a unique (up to
sign) $\frac {2\pi}k$-periodic solution. Moreover, all the
solutions are sign changing.
\end{proposition}

\subsubsection{Assume $\mu<0$}

This case was studied by Bidaut-V\'eron and Bouhar \cite{BVB}.
They obtained the following

\begin{proposition}\label{pro3} (\cite[Lemma 1.1 and
Lemma 1.2]{BVB}) Assume that $\mu<0$, $q>1$ and $\varepsilon=-1$.
Then
\begin{enumerate}
\item For all positive integer $k$ there exists a unique (up to sign) sign
changing $\frac {2\pi}k$-periodic solution of the ODE
(\ref{ode2}).

\item The ODE (\ref{ode2}) admits positive $2\pi$-periodic
solutions if and only if $-\mu (q-1)>1$. Moreover, in this case,
for any integer $1< k< \sqrt{-\mu(q-1)}$, the ODE (\ref{ode2})
admits a unique positive $\frac {2\pi}k$-periodic solution.
\end{enumerate}
\end{proposition}

\subsection{Case $q>1$ and $\varepsilon=+1$}

\subsubsection{Assume $\mu\le0$}\label{311}
This is an obvious case since the Laplace operator $w\mapsto w''$ is nonpositive on the circle.
Therefore, the only periodic solution of (\ref{ode2}) is the trivial one $w=0$.

\subsubsection{Assume $\mu>0$} This case
was studied by Chafee and Infante \cite{CI}. Here we give a
different approach based on the analysis of the dynamical system
$(S)$.

A direct calculation shows that the critical points of
$(\mathcal{S})$ are the origin and the two points $(-c,0)$ and
$(c,0)$, where $c=\mu^{\frac 1{q-1}}$. The origin is a center
while the two others are saddle points.

A classical study of the dynamical system gives the following

\begin{lemma}\label{pro4} Assume that $\mu>0$,
$q>1$ and $\varepsilon=+1$. Then the dynamical system $(S)$ satisfies the following properties:
\begin{enumerate}
\item There exists a unique heteroclinic orbit emanating from
$(-c,0)$ which tends to $(c,0)$ as $t\to \infty$ in the upper half
plane and one heteroclinic orbit emanating from $(c,0)$ which tends
to $(-c,0)$ as $t\to \infty$ in the lower half plane. The
equations of these orbits are given by
$$y^2=-\mu \left( x^2-c^2\right)+\frac2{q+1}\left[|x|^{q+1}-c^{q+1} \right].$$

\item Every point in the open bounded region delimited by these
two heteroclinic orbits belongs to a periodic orbit which turns
around the origin.

\item The period function $y_0\in(0,\gamma ) \mapsto T(y_0)$,
where $T(y_0)$ is the period of the solution passing through the
point $(0,y_0)$ and $\gamma=\sqrt{\mu c^2-\frac2{q+1}c^{q+1}}$, is
increasing with $T(0,\gamma )=\left(\frac {2\pi}{\sqrt
\mu},\infty\right)$.
\end{enumerate}
\end{lemma}

\noindent\textbf{Proof.} By classical arguments one can show that for
 $y_0>0$ small enough, orbits emanating from $(0,y_0)$ turn around the origin, while
for $y_0>0$ large enough, they still contained in a half plane $\left\{y>\alpha\right\}$
for some $\alpha>0$. This ensures the existence, and by
regularity, the uniqueness, of the heteroclinic orbits described
above.

To study the period function, one follows the same arguments as in the proof of Lemma \ref{pro1}
 (in this case, the expression of $U$ must be $U(x)=\mu x^2-\frac 2{q+1}x^{q+1}$).
We get
\begin{eqnarray*}
T(s)&=&4\int_0^{x(s)}\frac{dx}{\sqrt{s^2-U(x)}}=4\int_0^1\frac{x(s)d\tau}{\sqrt{U(x(s))-U(x(s)\tau)}}\\
&=& 4\int_0^1\frac{d\tau}{\sqrt{\mu (1-\tau^2)-\frac
2{q+1} x(s)^{q-1}(1-\tau^{q+1})}}.\\
\end{eqnarray*}
Hence, $\lim_{s\to 0}T(s)=
\frac 4{\sqrt\mu}\int_0^1\frac{d\tau}{\sqrt{1-\tau^2}}=\frac
{2\pi}{\sqrt\mu}$. Since $x(s)\to c=\mu^{\frac 1{q-1}}$ as $s\to \gamma= \sqrt{U(c)}$ (recall that $U(x(s))=s^2$)
and that the integral
$$\int_0^1\frac{d\tau}{\sqrt{\mu(1-\tau^2)-\frac{2\mu}{q+1}(1-\tau^{q+1})}}$$
is infinite, we deduce that $\lim_{s\to \gamma}T(s)=\infty$. \qeda

\bigskip

As a consequence we have

\begin{proposition}\label{pro5}
Assume that $q>1$, $\mu>0$ and $\varepsilon=+1$. The ODE
(\ref{ode2}) admits $2\pi$-periodic solutions if and only if
$\mu>1$. Moreover, in this case, for all positive integer $k<
\sqrt{\mu}$, the ODE (\ref{ode2}) admits a unique (up to sign)
$\frac{2\pi}k$-periodic solution which is sign changing.
\end{proposition}

\subsection{Case $0<q<1$ and $\varepsilon=-1$}

\subsubsection{Assume $\mu\ge0$}

The system $(S)$ is of class $C^1$ in
$\mathbb{R}^2\backslash\{O\}$, the origin being its unique
singular point. However the analysis is almost the same as for
$q>1$. Indeed, using Bendixson-Poincar\'e theory, closed orbits
must turn around the origin.

Using arguments like in the proof of Lemma \ref{pro1} and
the Lyapounov function
$$E(w,w')=\mu w^2+\frac 2{q+1}w^{q+1}+w'^2,$$
one can show the following

\begin{lemma}\label{pro6}
The period function $y_0\in (0,\infty) \mapsto T(y_0)$, where
$T(y_0)$ is the period of the solution of
$(\mathcal{S})$ passing through the point $(0,y_0)$, is increasing with
$$T(0,\infty)=\left\{\begin{array}{ll} (0, \frac{2\pi}{\sqrt \mu}) & \mbox{if $\mu >0$}\\
(0,\infty)&\mbox{if $\mu = 0$.}
\end{array} \right.$$
\end{lemma}

Consequently, one obtains the following

\begin{proposition}\label{pro7}
Assume that $0<q<1$, $\mu\ge0$ and $\varepsilon=-1$. Then for all
integers $k>\sqrt\mu$, the ODE (\ref{ode2}) admits a unique (up to
sign) $\frac {2\pi}k$-periodic solution. Moreover, all these
solutions are sign changing.
\end{proposition}

\subsubsection{Assume $\mu<0$}\label{323}

A direct calculation shows that the critical points of
$(\mathcal{S})$ are the two saddle points $(-c,0)$ and $(c,0)$, where
$c=(-\mu)^{\frac 1{q-1}}$. The origin is a singular point.

The proof of the following lemma is a slight modification of
that of Lemma \ref{pro4}.

\begin{lemma}\label{pro8} Assume that $\mu<0$, $0<q<1$ and
$\varepsilon=-1$. Then the dynamical system $(S)$ satisfies the following properties:
\begin{enumerate}
\item There exists a unique heteroclinic orbit emanating from
$(-c,0)$ which tends to $(c,0)$ as $t\to \infty$ in the upper half
plane and one heteroclinic orbit emanating from $(c,0)$ which tends
to $(-c,0)$ as $t\to \infty$ in the lower half plane. The
equations of these orbits are given by
$$y^2=-\mu \left( x^2-c^2\right)+\frac2{q+1}\left[x|x|^q-c^{q+1} \right].$$

\item Every point, except the origin, in the bounded region delimited by
these two heteroclinic orbits belongs to a periodic orbit which
turns around the origin.

\item The period function $y_0\in(0,\gamma)\mapsto T(y_0)$, where
$T(y_0)$ is the period of the solution passing through the point
$(0,y_0)$ and $\gamma=\sqrt{\mu c^2+\frac 2{q+1} c^{q+1}}$, is
increasing with
$T\left(0,\gamma\right)=\left(0,\infty\right)$.
\end{enumerate}
\end{lemma}

\noindent\textbf{Proof.} The existence and uniqueness of the heteroclinic orbits
rely on the same observation as in the proof of Lemma \ref{pro4}. Moreover,
using the same arguments as in the proof of Lemma \ref{pro1}, one obtains for
the period function (here, we have $U(x)=\mu x^2+\frac 2{q+1} x^{q+1}$),
\begin{eqnarray*}
T(s)&=&4\int_0^{x(s)}\frac{dx}{\sqrt{s^2-U(x)}}=4\int_0^1\frac{x(s)d\tau}{\sqrt{U(x(s))-U(x(s)\tau)}}\\
&=& 4x(s)^{\frac{1-q}2}\int_0^1\frac{d\tau}{\sqrt{\mu
x(s)^{1-q}(1-\tau^2)+\frac 2{q+1}(1-\tau^{q+1})}}.
\end{eqnarray*}
As before, we deduce that $\lim_{s\to 0}T(s)= 0$ and, since
$x(s)\to c= (-\mu)^{\frac 1{q-1}}$ as $s\to \gamma=\sqrt{U(c)}$,
$$\lim_{s\to \gamma}T(s)= \frac{4}{\sqrt{-\mu}}
\int_0^1\frac{d\tau}{\sqrt{\frac 2{q+1}(1-\tau^{q+1})-(1-\tau^2)}}=+\infty.$$ \qeda

\bigskip

As a consequence we have

\begin{proposition}\label{pro9}
Assume that $\mu<0$, $0<q<1$ and $\varepsilon=-1$. For all
positive integer $k$, the ODE (\ref{ode2}) admits a unique (up to sign)
$\frac{2\pi}k$-periodic solution. Moreover, all these solutions
are sign changing.
\end{proposition}

\subsection{Case $0<q<1$ and $\varepsilon=+1$}

\subsubsection{Assume $\mu\le0$}

As in subsection \ref{311}, in this case the only periodic solution is the
trivial one.

\subsubsection{Assume $\mu>0$}

The critical points of the system $(\mathcal{S})$
 are the two points $(-c,0)$ and $(c,0)$, with $\displaystyle c=\mu^{\frac 1{q-1}}$, and they
are centers.

\begin{remark}\label{r1}
Notice that $Q'(x)=-\mu+ q |x|^{q-1}$ is negative for $x>c$.
Hence, $Q(x)$ is positive for $0<x<c$ and negative
and decreasing for $x>c$.
\end{remark}

\begin{lemma}\label{l3}
The orbits of the dynamical system $(S)$ never admit horizontal or
vertical asymptotes.
\end{lemma}
\noindent{\bf Proof.} The existence of a vertical asymptote for
$(x(t), y(t))$ means that there exists $a\in {\mathbb R}$ such
that $x(t)\to a$ and $y(t) \to \pm \infty$ as $t\to \infty$, which
contradicts the equation $x'(t)=y(t)$.

In a similar manner, we can exclude horizontal asymptotes by noticing that
$y'(t)=Q(x(t))$ with $Q(x)\to \pm\infty$ as $x\to\pm\infty$. \qeda

\begin{lemma}\label{pro10}

\begin{enumerate}

\item Every orbit which intersects the $y$-axis at $y_0\ne 0$ is
periodic and turns around the origin and the two critical points.

\item Every orbit which intersects the $x$-axis at $x_0\in(-c,c)$,
$x_0\ne 0$ is periodic and turns around the critical point $(c,0)$
if $x_0>0$, and around the point $(-c,0)$ if $x_0<0$.

\end{enumerate}
\end{lemma}

\noindent{\bf Proof. 1.} Using time shift invariance and symmetry,
one can assume that the orbit starts at $(0,y_0)$ with $y_0> 0$.
From $y'=Q(x)$ and Lemma \ref{l3}, one deduces that $y$ increases
for $0<x<c$, then decreases for $x>c$ and intersects the $x$-axis
at finite time. Applying
Lemma \ref{l2} the orbit is periodic.\\
{\bf 2.} By the same analysis, every orbit which starts at $(x_0,0)$
with $0<x_0<c$ must intersect the $x$-axis at finite time and, hence, is
periodic. By symmetry the same result holds for $-c<x_0<c$.\qeda

\begin{lemma}\label{pro11}
There exists exactly one homoclinic-like orbit emanating from the
origin in the right (resp. left) half-plane and enclosing the critical
point $(c,0)$ (resp. $(-c,0)$). Moreover the equation of
this orbit is given by
\begin{equation}\label{homo1}y^2+\mu x^2=\frac 2{1+q}|x|^{1+q}.\end{equation}
\end{lemma}

\begin{remark}
In the particular case $\mu=\frac4{(1-q)^2}$, which corresponds to
the Euclidean metric, the orbit (\ref{homo1}) could be
parametrized by $(x(t),x'(t))$ where
\begin{equation}\label{homo} x(t):=x_*|\sin
(t)|^{\sqrt\mu} \end{equation} with $\displaystyle
x_*:={\left(\frac 2{(1+q)\mu} \right)}^{\frac 1{1-q}}$.

\end{remark}

\noindent{\bf Proof.} \underline{Existence}. Integrating the
differential equation
$$\frac{dy}{dx}=\frac{-\mu x+x^q}y$$
from $0$ to $t$, we get
\begin{equation}\label{eq}
y^2+U(x)=y^2(0)+U(x(0)),
\end{equation}
where $U(x):=\mu x^2-\frac 2{q+1}x^{q+1}$. Thus the equation of
the orbit passing through the origin is given by
$$y^2=-U(x)=-\mu x^2+\frac 2{q+1}x^{q+1}.$$
\noindent\underline{Uniqueness}. Since the origin is a singular point, we
cannot use classical theorems to deduce uniqueness of the special
homoclinic orbits given by (\ref{homo1}). To show uniqueness,
consider the change of variables $X:=\varphi(x)$ and $Y:=y$, where
$y=\varphi(x)$ is the equation of the upper part of the positive
special homoclinic orbit defined by (\ref{homo1}). Then $X'=\frac
d{dt}\left (\varphi(x)\right) =x'\varphi'(x)=y\varphi'(x)$ with
$\varphi'(x) = Q(x)/\varphi(x)$. On the other hand,
$Y'=y'=Q(x)$. Thus $X'=Q(x)y/\varphi(x)=Y'Y/X$, i.e. $\frac
d{dt}\left (X^2\right)=2X'X=2Y'Y=\frac d{dt}\left (Y^2\right)$.
Therefore
$X^2(t)-Y^2(t)=X^2(0)-Y^2(0)$ which gives the uniqueness. \qeda

\medskip

\noindent\textbf{Variation of the period}

\medskip

\noindent{\it Sign changing solutions}
\begin{lemma}\label{pro12}
The period function $y_0\in\left(0,\infty\right)\mapsto T(y_0)$,
where $T(y_0)$ is the period of the solution passing through the
point $(0,y_0)$, is decreasing and
$T\left(0,\infty\right)=\left(
\frac{2\pi}{\sqrt\mu},\frac{4\pi}{(1-q)\sqrt\mu}\right)$.
\end{lemma}
\noindent{\bf Proof.} We follow the same arguments as in the proof
of Lemma \ref{pro1}, we obtain, with the same notations. In this
case, the expression of $U$ is $U(x)=\mu x^2-\frac 2{q+1}x^{q+1}$.
We get
\begin{eqnarray*}
T(s)&=&4\int_0^{x(s)}\frac{dx}{\sqrt{s^2-U(x)}}=4\int_0^1\frac{x(s)d\tau}{\sqrt{U(x(s))-U(x(s)\tau)}}\\
&=& 4\int_0^1\frac{d\tau}{\sqrt{\mu (1-\tau^2)-\frac
2{q+1} x(s)^{q-1}(1-\tau^{q+1})}}.
\end{eqnarray*}
Since $x(s)\to\infty$ as $s\to\infty$, we have
$$\lim_{s\to \infty}T(s)=\frac 4{\sqrt\mu}\int_0^1\frac{d\tau}{\sqrt{1-\tau^2}}=\frac
{2\pi}{\sqrt\mu}.$$ On the other hand, when $s$ goes to $0$,
$x(s)$ tends to $x_*$, the positive solution of $U(x_*)=0$ (i.e.
$x_*^{q-1}=\frac{q+1}2\mu$). Therefore,
$$\lim_{s\to 0}T(s)=4\int_0^{x_*}\frac{dx}{\sqrt{-U(x)}},$$
and, using the change of variable $x=x_*\sin
^{\frac2{1-q}}t$, we get
$$\lim_{s\to 0}T(s)=\frac{8}{1-q}\int_0^{\frac{\pi}{2}}
\frac{\cos t\, dt}{\sqrt{\frac{2}{q+1}x_*^{q-1}-\mu \sin^2t}}=
\frac{8}{(1-q)\sqrt{\mu}}\int_0^{\frac{\pi}{2}} dt=\frac{4\pi}{(1-q)\sqrt{\mu}}.$$
\hfill\qeda


\newpage

\noindent{\it Positive solutions}

Observe that the homoclinic-like solution given by (\ref{homo}) is
$\frac{2\pi}{(1-q)\sqrt{\mu}}$-periodic. On the other hand, the
linearized equation at $(c,0)$ is
$$w''+\mu(1-q)w=0.$$
Roughly speaking the period function takes the value
$\frac{2\pi}{\sqrt{\mu(1-q)}}$ at the point orbit $(c,0)$ and
$\frac{2\pi}{(1-q)\sqrt{\mu}}$ at the homoclinic-like orbit.

\begin{lemma}\label{pro13}
The period function $s\in(0,\gamma)\mapsto T(s)$, where $T(s)$ is
the period of the solution passing through the point $(c,s)$ and
$\gamma>0$ such that $\gamma^2=-\mu c^2+\frac 2{q+1}c^{q+1}$, is
increasing with
$T\left(0,\gamma\right)=\left(\frac{2\pi}{\sqrt{\mu(1-q)}},\frac{2\pi}{(1-q)\sqrt{\mu}}\right)$.
\end{lemma}
\noindent\textbf{Proof.} For any $s\in (0,\gamma)$ we denote by $(y(s),0)$ and $(z(s),0)$,
with $0<y(s)<c<z(s)<x_*$, the intersection points of
the orbit passing through the point $(c,s)$ with the $x$-axis.
Following an idea of \cite{ACP} (see also \cite{BVB}), since
$y^2=s^2-U(x)+U(c)$, we have
\begin{eqnarray*}T(s)/2&=&\int_c^{z(s)}\frac{du}{\sqrt{s^2-U(u)+U(c)}}
-\int_c^{y(s)}\frac{du}{\sqrt{s^2-U(u)+U(c)}}\\
&=&\int_0^1(z'(st)-y'(st))\frac{dt}{\sqrt{1-t^2}}.\end{eqnarray*}
In particular, since $\gamma^2=-U(c)$,$$\lim_{s\to \gamma} T(s) =
\int_0^{x_*}\frac{2du}{\sqrt{\gamma^2-U(u)+U(c)}}=
\int_0^{x_*}\frac{2du}{\sqrt{-U(u)}}=\frac{2\pi}{(1-q)\sqrt{\mu}}$$
(see the end of the proof of Lemma \ref{pro12}). On the
other hand, since $y(s)$ and $z(s)$ satisfy
$$0=U(y(s))-U(c)-s^2=U(z(s))-U(c)-s^2,$$
we easily get $$\lim_{s\to 0}\left(z'(s)-y'(s)\right)=\frac2
{\sqrt{\mu(1-q)}}.$$ Therefore $$\lim_{s\to 0}T(s)=
\frac{2\pi}{\sqrt{\mu(1-q)}}.$$

Now, to show the monotony of the function $T$, we first observe that
$$
\frac d{ds}T(s)=2\int_0^1\left(z''(st)-y''(st)
\right)\frac{tdt}{\sqrt{1-t^2}}.
$$
Thus, it suffices to show that $z''-y''$ is nonnegative. From the
equation $U(z(s))-U(c)-s^2=0$ we get $z'(s)=2s/U'(z(s))$ and, hence,
$$z''(s)=2\frac{U'^2(z(s))+2(U(c)-U(z(s)))U''(z(s))}{U'^3(z(s))}.$$
Since $z(s)\to c$ as $s\to 0$, $U$ is smooth, $U'(c)=0$ and
$U''(c)\ne0$, expanding $U$ up to order 4 in the neighborhood of
$c$ shows that
$$\lim_{s\to 0} 2\frac{U'^2(z(s))+2(U(c)-U(z(s)))U''(z(s))}{U'^3(z(s))}= -
\frac{2U^{(3)}(c)}{U''^2(c)}.$$
Thus,
$$\lim_{s\to 0}z''(0)=-\frac{2U^{(3)}(c)}{3U''^2(c)}=\lim_{s\to 0}y''(0).$$
Let us show that \begin{equation}\label{eqzy}z''(s)\ge-
\frac{2U^{(3)}(c)}{U''^2(c)}\ge y''(s).\end{equation}
For the first inequality, it suffices to show that
$$F:=U'^2+2U''[U(c)-U(z)]+\frac{U'^3U^{(3)}(c)}{3U''^2(c)}\ge 0$$
on $(c,x_*)$. In fact we will show that this is true even on
$(0,+\infty)$. Indeed, notice first that
$$F'=2U^{(3)}[U(c)-U]+\frac{U^{(3)}(c)}{U''^2(c)}U'^2U''\le 0$$
on $(0,q^{\frac 1{1-q}}c)$. On the other hand, set $H:=-\frac F{U''}$ on
$(q^{\frac 1{1-q}}c, \infty)$, then we have $H'=-K\frac
{U'^2U^{(3)}}{U''^2}$, where
$$K=-1+\frac{U^{(3)}(c)}{3U''^2(c)}\left[3\frac {U''^2}{U^{(3)}}
-U'\right].$$
Now,
$$K'=\frac{U^{(3)}(c)U''}{3U''^2(c)(U^{(3)})^2}\left[5(U^{(3)})^2-3U^{(4)}U''\right]>0$$
with $K(c)=0$, then $K$ is nonpositive on $(q^{\frac 1{1-q}}c,c)$
and nonnegative on $(c,\infty)$. Hence, $H'$ is nonnegative on
$(q^{\frac 1{1-q}}c,c)$ and nonpositive on $(c,\infty)$ and since
$H(c)=0$ we get $H\le 0$. Thus $F\ge 0$ on $(q^{\frac 1{1-q}}c,
\infty)$, Therefore, $F\ge0$ on $(0,\infty)$. In a similar manner,
one can show the second inequality in (\ref{eqzy}). \qeda

The following proposition summarizes the consequences of Lemmas
\ref{pro11}, \ref{pro12} and \ref{pro13}.

\begin{proposition}\label{pro14} Assume that $\mu>0$,
$\varepsilon=+1$ and $0<q<1$. Then

\begin{enumerate}
\item The ODE (\ref{ode2}) admits nonconstant $2\pi$-periodic sign
changing solutions if and only if $\mu>1$. Moreover, in this case,
for all integer $k\in\left(\frac {1-q}2\sqrt\mu,\sqrt\mu\right)$ the ODE
(\ref{ode2}) admits a unique (up to sign) sign changing $\frac {2\pi}
k$-periodic solution.

\item The ODE (\ref{ode2}) admits nonconstant $2\pi$-periodic
positive solutions if and only if $\left((1-q)\sqrt\mu,
\sqrt{\mu(1-q)}\right)\cap \mathbb{N}\ne\emptyset$. Moreover, in this
case, for all integer $(1-q)\sqrt\mu< k< \sqrt{\mu(1-q)}$, the ODE
(\ref{ode2}) admits a unique positive $\frac{2\pi}k$-periodic
solution.

\item There exists a unique
$\frac{2\pi}{(1-q)\sqrt{\mu}}$-periodic nonnegative solution of
(\ref{ode2}). It vanishes only once in a period and its orbit is
given by (\ref{homo1}).

\end{enumerate}
\end{proposition}


\section{Applications}
In this section we apply the results above to the
existence problem of non-radial pseudo-radial solutions of the PDE
(\ref{edp}) and its particular case (\ref{conf}).

\subsection{Case of a metric $g=a^2(r)g_{{\mathbb R}^2}$}

In this case $b(r)=ra(r)$. Note that
$\Delta^gu=\frac1{a^2(r)}\Delta u$, where $\Delta$ is the
Euclidean Laplacian. Hence, the PDE
$$\Delta^g u=\varepsilon u|u|^{q-1}$$
is equivalent to the following one
\begin{equation}\label{eDa}
\Delta u=\varepsilon a^2(|x|)u|u|^{q-1}
\end{equation}
that we consider on a disc or an annular domain
$\Omega:=\{0\le R_1<|x|<R_2\le +\infty\}$ of $\mathbb{R}^2$.

\begin{theorem}\label{th2}
If the equation (\ref{eDa}) admits a non-radial pseudo-radial
solution in $\Omega$, then the function $r\mapsto a(r)$ has one of
the three following forms:
 $$a(r)=r^{-1}{\left[
Mr^\alpha+Nr^{-\alpha}\right]}^{\frac{1-q}2},$$ where $\alpha>0$
and $M,N\in\mathbb{R}$ are such that $Mr^\alpha+Nr^{-\alpha}>0$ on
$(R_1,R_2)$,
$$a(r)=r^{-1}{\left[ M+N\ln r \right]}^{\frac{1-q}2},$$
where $M,N\in\mathbb{R}$ are such that $M+N\ln r>0$ on
$(R_1,R_2)$, or
$$a(r)=r^{-1}{\left[ M\cos(\alpha\ln r) +
N\sin(\alpha\ln r) \right]}^{\frac{1-q}2},
$$ where $\alpha>0$ and $M,N\in\mathbb{R}$.

\end{theorem}
\noindent\textbf{Proof.} Applying Theorem \ref{main}, the equation
(\ref{eq2}) becomes (with $c=b/a=r$ and $h=b^{\frac{2}{1-q}}$)
\begin{equation}\label{e2bis}{(rh')}'=rh''+h'=\mu \frac h r.\end{equation}
The general solution of (\ref{e2bis}) is given by
$$h(r)=\left\{\begin{array}{ll}
Mr^{\sqrt\mu}+Nr^{-\sqrt{\mu}} &\mbox{if }\mu>0,\\
\\
M+N\ln r&\mbox{if }\mu=0,\\
\\
M\cos(\sqrt{-\mu}\ln r) + N\sin(\sqrt{-\mu}\ln r) &\mbox{if
}\mu<0.
\end{array} \right.
$$
 \qeda

Hence, for our purpose, the only relevant equations of the form
$\Delta u=\varepsilon a^2(|x|)u|u|^{q-1}$ in $\Omega$ are
\begin{equation}\label{eDa1}
\Delta u=\varepsilon |x|^{-2}{\left[
M|x|^\alpha+N|x|^{-\alpha}\right]}^{1-q}u|u|^{q-1},
\end{equation}
\begin{equation}\label{eDa2}
\Delta u=\varepsilon |x|^{-2}{\left[ M+N\ln
|x|\right]}^{1-q}u|u|^{q-1},
\end{equation}
and
\begin{equation}\label{eDa3}
\Delta u=\varepsilon |x|^{-2}{\left[ M\cos(\alpha\ln |x|) +
N\sin(\alpha\ln |x|)\right]}^{1-q}u|u|^{q-1},
\end{equation}
where $\alpha$, $M$ and $N$ are as in Theorem \ref{th2}.

For (\ref{eDa1}), the set of pseudo-radial non-radial solutions
can be described in the following way:

\begin{theorem}\label{th3}
\begin{enumerate}
\item If $\varepsilon=+1$ and $0<\alpha\le1$, then any
pseudo-radial solution of (\ref{eDa1}) is radial.

\item Define, for $\varepsilon=-1$ or $\alpha>1$, the subset
$\mathcal{A}(\varepsilon,q,\alpha)\subset {\mathbb N} $ by
\begin{eqnarray*}
\mathcal{A}(\varepsilon,q,\alpha)&=&\left\{
 \begin{array}{ll}
 (\alpha,+\infty)\cap\mathbb{N},&\mbox{if }\varepsilon=-1\\
 (0,+\infty)\cap\mathbb{N}&\mbox{if }\varepsilon=+1, q>1\mbox{ and }\alpha>1\\
 (\alpha(1-q)/2,\alpha)\cap\mathbb{N}&\mbox{if }\varepsilon=+1, 0<q<1\mbox{ and
 }\alpha>1,
 \end{array} \right.
\end{eqnarray*}
Then the set of sign changing non-radial pseudo-radial solutions of (\ref{eDa1}) is
parameterized by $\mathcal{A}(\varepsilon,q,\alpha)$ in the sense
that it consists in the functions
$$u_k(x)=\pm [M|x|^\alpha+N|x|^{-\alpha}]w_k(\theta),\; \; k\in \mathcal{A}(\varepsilon,q,\alpha),$$
where, for every $ k\in \mathcal{A}(\varepsilon,q,\alpha)$, $w_k$
is the unique (up to sign) sign changing $\frac{2\pi}k$-periodic
solution of (\ref{ode2}). \item Equation (\ref{eDa1}) admits
positive non-radial pseudo-radial solutions if and only if
$\varepsilon = +1$, $0<q<1$ and $\mathcal{B}(+1,q,\alpha):=
(\alpha(1-q),\alpha\sqrt{1-q})\cap\mathbb{N}\neq \emptyset$.
Moreover, all these solutions are of the form
$$u_k(x)= [M|x|^\alpha+N|x|^{-\alpha}]v_k(\theta),$$
where, for every $ k\in \mathcal{B}(+1,q,\alpha)$, $v_k$ is the
unique positive $\frac{2\pi}k$-periodic solution of (\ref{ode2}).
\item If $(1-q)\alpha\in {\mathbb N} $, then equation (\ref{eDa1})
admits a unique nonnegative pseudo-radial solution given by
$$ u(x)= [M|x|^\alpha+N|x|^{-\alpha}]v(\theta),$$
where $v$ is the unique nonnegative periodic solution of
(\ref{ode2}) (corresponding to the homoclinic-like orbit) which
vanishes once in a period.
\end{enumerate}
\end{theorem}

Concerning (\ref{eDa2}), we have the following

\begin{theorem}\label{th4}
\begin{enumerate}
\item If $\varepsilon=+1$, then, for every $ q>0$, $q\ne1$, any
pseudo-radial solution of (\ref{eDa2}) is radial.

\item If $\varepsilon=-1$, then, for every $ q>0$, $q\ne 1$, the
set of non-radial pseudo-radial solutions of (\ref{eDa2}) is
parameterized by ${\mathbb N} ^*$ and consists in the functions
$$u_k(x)=\pm [M+N\ln |x|]w_k(\theta),$$
where, for every $ k\in\mathbb{N}^*$, $w_k$ is the unique (up to
sign) sign changing $\frac{2\pi}k$-periodic solution of
(\ref{ode2}).
\end{enumerate}
\end{theorem}

\begin{remark}
A particular case of (\ref{eDa1}) is
\begin{equation}\label{eDap}
\Delta u =\varepsilon |x|^p u|u|^{q-1}, \quad p\in{\mathbb R}\backslash\left\{-2\right\}.
\end{equation}
One can apply Theorem \ref{th3} with
$\alpha=\Big|\frac{p+2}{q-1}\Big|$, $(M,N)=(1,0)$ if
$\frac{p+2}{q-1}<0$ and $(M,N)=(0,1)$ if $\frac{p+2}{q-1}>0$.

The case $p=-2$ corresponds to (\ref{eDa2}) with $M=1$ and $N=1$ and is covered by Theorem \ref{th4}.
\end{remark}

Concerning (\ref{eDa3}), we have the following

\begin{theorem}\label{th5}
\begin{enumerate}
\item If $\varepsilon=+1$, then, for every $ q>0$, $q\ne1$, any
pseudo-radial solution of (\ref{eDa3}) is radial.

\item For $\varepsilon=-1$, the set of sign changing non-radial
pseudo-radial solutions of (\ref{eDa3}) is parameterized by
${\mathbb N} ^*$ and consists in the functions
$$u_k(x)=\pm [M\cos(\alpha\ln |x|) +
N\sin(\alpha\ln |x|) ]w_k(\theta),$$ where, for every
$k\in\mathbb{N}^*$, $w_k$ is the unique (up to sign) sign changing
$\frac{2\pi}k$-periodic solution of (\ref{ode2}).

\item Equation (\ref{eDa3}) admits positive non-radial
pseudo-radial solutions if and only if $\varepsilon = -1$, $q>1$
and $\alpha\sqrt{q-1}>2$. Moreover, all these solutions are of the
form
$$u_k(x)= [M\cos(\alpha\ln |x|) +
N\sin(\alpha\ln |x|) ]v_k(\theta),$$ where, for every $k\in
(1,\alpha\sqrt{q-1})\cap \mathbb{N}$, $v_k$ is the unique positive
$\frac{2\pi}k$-periodic solution of (\ref{ode2}).
\end{enumerate}
\end{theorem}

\subsection{The spherical case: $a(r)=1$ and $b(r)=\sin r$}
In this case, condition (\ref{conditiong}) implies
\begin{equation}\label{eqbb}
\mu\frac {1-q}2=\frac2{1-q}\cos^2r-\sin^2 r,
\end{equation}
which is only possible if $q=3$ and $\mu=1$. Applying Propositions
\ref{pro2} and \ref{pro5} we get

\begin{theorem}\label{th6}
Let $\Omega\subset{\mathbb S}^2$ be a rotationally symmetric
domain of the standard sphere and consider in $\Omega$ the following equation:
\begin{equation}\label{edpS2}
\Delta_{{\mathbb S}^2} u=\varepsilon u|u|^{q-1},
\end{equation}
where $\Delta_{{\mathbb S}^2} $ is the standard Laplacian of ${\mathbb S}^2$.
\begin{itemize}
\item[i)] If $\varepsilon=+1$ or $q\ne3$, then any pseudo-radial solution
of (\ref{edpS2}) is radial.
\item[ii)] The equation $\Delta_{{\mathbb S}^2} u=- u^3 $ admits
infinitely many non-radial pseudo-radial solutions in $\Omega$ which are all of the form
$$u_k(x)=\frac{w_k(\theta)}{\sin r},$$
where, for any integer $ k\ge 2$, $w_k$ is the unique (up to sign)
$\frac{2\pi}k$-periodic solution of (\ref{ode2}).
\end{itemize}

\end{theorem}

\subsection{The hyperbolic metric case: $a(r)=1$ and $b(r)=\sinh r$}
Here also, condition (\ref{conditiong}) is satisfied if and only
if $q=3$ and $\mu=1$. Like in the spherical case, we get the
following

\begin{theorem}\label{th7} Let $\Omega\subset{\mathbb H}^2$ be a rotationally symmetric
domain of the hyperbolic plane and consider in $\Omega$ the following equation:
\begin{equation}\label{edpH2}
\Delta_{{\mathbb H}^2} u=\varepsilon u|u|^{q-1},
\end{equation}
where $\Delta_{{\mathbb H}^2} $ is the Laplacian of ${\mathbb H}^2$.
\begin{itemize}
\item[i)] If $\varepsilon=+1$ or $q\ne3$, then any pseudo radial solution
of (\ref{edpH2}) is radial.
\item[ii)] The equation $\Delta_{{\mathbb H}^2} u=- u^3 $ admits
infinitely many non-radial pseudo-radial solutions in $\Omega$ which are all of the form
$$u_k(x)=\frac{w_k(\theta)}{\sinh r},$$
where, for any integer $ k\ge 2$, $w_k$ is the unique (up to sign)
$\frac{2\pi}k$-periodic solution of (\ref{ode2}).
\end{itemize}\end{theorem}

\subsection{Case of a metric conformal to the cylindrical one: $a=b$}
The standard metric of the cylinder $C=(R_1,R_2)\times{\mathbb
S}^1$ is the product metric $dr^2+d\theta^2$, its Laplacian is
given by $\Delta_C=\frac{\partial^2}{\partial
r^2}+\frac{\partial^2}{\partial \theta^2}$. For a conformal metric
$g=a^2(r)[dr^2+d\theta^2]$, the associated Laplacian is
$\Delta^g=a^{-2}(r)\Delta_C$. The equation
\begin{equation}\label{eDa4}
\Delta^S u=\varepsilon a(r)^{2}u|u|^{q-1},
\end{equation}
is then equivalent to
$$\Delta^gu=\varepsilon u|u|^{q-1}.$$
Condition (\ref{conditiong}) gives, with $a=b$,
$${\left[a^{\frac{1+q}{1-q}}a'\right]}'=\frac{\mu}2(1-q)a^{\frac2{1-q}}.$$
Setting $\gamma=\frac{1+q}{1-q}$, this last equation becomes
$$(a^{\gamma+1})''=\mu a^{\gamma+1}$$
which gives
$$a^{\gamma+1}(r)=\left\{\begin{array}{ll}
A\cosh\sqrt\mu \,r+B\sinh\sqrt\mu \,r&\mbox{if }\mu>0\\
\\
A +Br&\mbox{if }\mu=0\\
\\
A\cos\sqrt{-\mu} \, r+B\sin\sqrt{-\mu}\, r&\mbox{if
}\mu<0,
\end{array}\right.$$
where $A$ and $B$ are such that the right hand side is positive on $(R_1,R_2)$.
Therefore, for our purpose, the only relevant equations of the form (\ref{eDa4})
are
\begin{equation}\label{eDa5}
\Delta_C u=\varepsilon {\left[
A\cosh\alpha r+B\sinh\alpha r\right]}^{1-q}u|u|^{q-1},
\end{equation}
\begin{equation}\label{eDa6}
\Delta u=\varepsilon {\left[ A+Br\right]}^{1-q}u|u|^{q-1},
\end{equation}
and
\begin{equation}\label{eDa7}
\Delta u=\varepsilon {\left[ A\cos\alpha r+B\sin\alpha r\right]}^{1-q}u|u|^{q-1},
\end{equation}
where $\alpha>0$.

The description of the set of non-radial pseudo-radial solutions
of these equations is the same as for equations (\ref{eDa1}),
(\ref{eDa2}) and (\ref{eDa3}) respectively. Indeed, the statement
of Theorem \ref{th3} remains valid for equation (\ref{eDa5})
provided that $M|x|^\alpha+N|x|^{-\alpha}$ is replaced by
$A\cosh\alpha r+B\sinh\alpha r$. Similarly, Theorem \ref{th4}
applies to (\ref{eDa6}) replacing $M+N\ln |x|$ by $A+Br$, and
Theorem \ref{th5} applies to (\ref{eDa7}) replacing
$M\cos(\alpha\ln |x|) + N\sin(\alpha\ln |x|)$ by $A\cos\alpha
r+B\sin\alpha r$.


\vskip 2cm

\begin{tabular}{ll}
Ahmad EL SOUFI &{}\hskip 1 cm{} Mustapha JAZAR\\
\\
Universit\'e de Tours &{}\hskip 1 cm{} Lebanese University\\
Laboratoire de Math\'ematiques&{}\hskip 1 cm{} Mathematics Department\\
et Physique Th\'eorique &{}\hskip 1 cm{} P.O.Box 155-012 \\
UMR 6083 du CNRS& {}\hskip 1 cm{} Beirut, Lebanon \\
Parc de Grandmont&\\
F- 37200 Tours, France&\\

\\
elsoufi@univ-tours.fr &{}\hskip 1 cm{} mjazar@ul.edu.lb\\

\end{tabular}

\end{document}